\documentclass[smallextended,envcountsect]{svjour3}
\smartqed
\usepackage{graphicx}
\journalname{JOTA}

 \usepackage{amssymb,amsmath}

 \newcommand{\Id}{\operatorname{Id}_n}
 \newcommand{\AC}{AC^\alpha}
 \newcommand{\rd}{\mathrm{d}}
 \newcommand{\diam}{\operatorname{diam}}
 \newcommand{\USet}{\mathbb{U}}
 \newcommand{\U}{\mathcal{U}}
 \newcommand{\Val}{\rho}
 \newcommand{\Image}{\mathcal{I}}
 \newcommand{\Vala}{\rho_\ast}
 \newcommand{\Valaeta}{\rho_\eta}
 \newcommand{\argmax}[1]{\underset{#1}{\operatorname{argmax} \,}}
 \newcommand{\argmin}[1]{\underset{#1}{\operatorname{argmin} \,}}
 \newcommand{\esssup}[1]{\underset{#1}{\operatorname{ess\,sup} \,}}

\begin{document}

\title{Optimal Control Problems with a Fixed Terminal Time in Linear Fractional-Order Systems}

\author{Mikhail Gomoyunov}

\institute{Mikhail Gomoyunov \at
    Krasovskii Institute of Mathematics and Mechanics of the Ural Branch of the Russian Academy of Sciences;
    Ural Federal University \\
    Ekaterinburg, Russia \\
    m.i.gomoyunov@gmail.com
}

\date{Received: 23.08.2019 / Accepted: date}

\maketitle

\begin{abstract}
    The paper deals with an optimal control problem in a dynamical system described by a linear differential equation with the Caputo fractional derivative.
    The goal of control is to minimize a Bolza-type cost functional, which consists of two terms: the first one evaluates the state of the system at a fixed terminal time, and the second one is an integral evaluation of the control on the whole time interval.
    In order to solve this problem, we propose to reduce it to some auxiliary optimal control problem in a dynamical system described by a first-order ordinary differential equation.
    The reduction is based on the representation formula for solutions to linear fractional differential equations and is performed by some linear transformation, which is called the informational image of a position of the original system and can be treated as a special prediction of a motion of this system at the terminal time.
    A connection between the original and auxiliary problems is established for both open-loop and feedback (closed-loop) controls.
    The results obtained in the paper are illustrated by examples.
\end{abstract}


\section{Introduction}

In the paper, we consider an optimal control problem in a dynamical system which motion is described by a linear differential equation with the Caputo fractional derivative of an order $\alpha \in (0, 1)$.
The goal of control is to minimize a given Bolza-type cost functional, which consists of two terms.
One of them evaluates the state vector of the system realized at a fixed terminal time $\vartheta$, and the other is an integral evaluation of a control on the whole time interval $[t_0, \vartheta]$.
We are interested in finding optimal, or, at least, $\varepsilon$-optimal, open-loop controls, as well as in constructing optimal feedback (closed-loop) controls, which are formalized within the framework of positional strategies \cite{Krasovskii_Subbotin_1988,Krasovskii_Krasovskii_1995} (see also \cite{Gomoyunov_2019_Trudy_Eng,Gomoyunov_2019_SIAM}).

In order to solve the problem, we propose an approach based on its reduction to some auxiliary optimal control problem in a dynamical system described by a first-order ordinary differential equation and further application of the methods and results of the optimal control theory widely developed for such systems.
The reduction relies on a suitable notion of a finite-dimensional informational image of an infinite-dimensional position \cite{Gomoyunov_2019_SIAM} of the original fractional-order system, which can be treated as a special prediction of a motion of this system at the time $\vartheta$.

This approach is closely related to a functional interpretation of control processes \cite{Krasovskii_1963} (see also \cite{Lukoyanov_Gomoyunov_2019_DGAA} and the references therein) and was previously developed for linear functional-differential systems of retarded \cite{Lukoyanov_Reshetova_1998_Eng,Gomoyunov_Plaksin_2015_IFAC} and neutral \cite{Gomoyunov_Lukoyanov_2018_Trudy_Eng} types and also for linear systems with control delays \cite{Gomoyunov_Lukoyanov_2012_Eng}.
However, in contrast to these studies, the auxiliary dynamical system obtained in the present paper may not satisfy the assumptions that are usually made in the optimal control theory.
More precisely, the right-hand side of the corresponding differential equation, in general, has a singularity at the time $\vartheta$ and, therefore, is unbounded.
This circumstance is explained by some special properties of the fundamental matrix solution of linear fractional-order differential equations (see, e.g., \cite{Idczak_Kamocki_2011,Bourdin_2018,Gomoyunov_2019_FCAA_2}).
In order to overcome this difficulty, we propose to introduce a small parameter $\eta \in (0, \vartheta - t_0)$ and shift the terminal time in the auxiliary optimal control problem from $\vartheta$ to $\vartheta_\eta = \vartheta - \eta$.
After such a modification, the auxiliary problem, on the one hand, can still be used to solve, at least approximately, the original problem, and, on the other hand, it already meets the typical assumptions from the optimal control theory.

The rest of the paper is organized as follows.
In Sect.~\ref{section_OCP}, we formulate the optimal control problem and study the question of finding optimal and $\varepsilon$-optimal open-loop controls.
We introduce the notions of a position of the system and its informational image, derive the auxiliary optimal control problem, and provide a connection between the original and auxiliary problems.
After that, we consider the auxiliary optimal control problem with the shifted terminal time $\vartheta_\eta$ and establish its connection with the original problem.
Sect.~\ref{section_OFC} is devoted to constructing optimal feedback controls.
We introduce the notion of a positional control strategy and show how to obtain an optimal positional control strategy on the basis of optimal positional strategies in the auxiliary problems with the shifted terminal times $\vartheta_\eta$.
Both sections contain illustrative examples.
The proofs of the statements are given in Appendix.

\section{Optimal Control Problem}
\label{section_OCP}

\subsection{Preliminaries}

Let $n \in \mathbb{N}$, $\alpha \in (0, 1)$, and $t_0$, $\vartheta \in \mathbb{R}$ such that $t_0 < \vartheta$ be fixed throughout the paper.
Let $\mathbb{R}^n$ and $\mathbb{R}^{n \times n}$ be the spaces of $n$-dimensional vectors and $(n \times n)$-matrices, and let $\Id \in \mathbb{R}^{n \times n}$ stand for the identity matrix.
By $\|\cdot\|$ and $\langle \cdot, \cdot \rangle$, we denote the Euclidean norm and the inner product in $\mathbb{R}^n$.
The corresponding norm in $\mathbb{R}^{n \times n}$ is also denoted by $\|\cdot\|$.

Let $t_\ast \in [t_0, \vartheta]$.
Let us consider a function $x: [t_0, t_\ast] \rightarrow \mathbb{R}^n$, for which we also use the notation $x(\cdot)$.
The (left-sided) Riemann--Liouville fractional integral $(I^\alpha x)(t)$ and Caputo fractional derivative $(^C D^\alpha x)(t)$ of the order $\alpha$ of $x(\cdot)$ at $t \in [t_0, t_\ast]$ are defined respectively by
\begin{equation*}
    \begin{array}{c}
        \displaystyle
        (I^\alpha x) (t)
        = \frac{1}{\Gamma(\alpha)} \int_{t_0}^{t} \frac{x(\tau)}{(t - \tau)^{1 - \alpha}} \, \rd \tau, \\[0.5em]
        \displaystyle
        (^C D^\alpha x)(t)
        = \frac{\rd}{\rd t} \big(I^{1 - \alpha} (x(\cdot) - x(t_0)) \big) (t)
        = \frac{1}{\Gamma(1 - \alpha)} \frac{\rd}{\rd t} \int_{t_0}^{t} \frac{x(\tau) - x(t_0)}{(t - \tau)^\alpha} \, \rd \tau,
    \end{array}
\end{equation*}
where $\Gamma$ is the gamma-function.
The basic properties of integrals and derivatives of fractional order can be found, e.g., in \cite{Samko_Kilbas_Marichev_1993,Kilbas_Srivastava_Trujillo_2006,Diethelm_2010}.
We say that $x(\cdot) \in \AC([t_0, t_\ast], \mathbb{R}^n)$, if there exists a (Lebesgue) measurable and essentially bounded function $\varphi: [t_0, t_\ast] \rightarrow \mathbb{R}^n$ such that $x(t) = x(t_0) + (I^\alpha \varphi) (t)$, $t \in [t_0, t_\ast]$.
Let us note that, in the case when $t_\ast = t_0$, the set $\AC([t_0, t_\ast], \mathbb{R}^n)$ can be identified with $\mathbb{R}^n$.

\subsection{Statement of the Problem}
\label{subsection_OCP_Statement_of_the_problem}

We consider a dynamical system which motion is described by the linear fractional differential equation
\begin{equation} \label{system}
    (^C D^\alpha x)(t) = A(t) x(t) + f(t, u(t)),
    \quad x(t) \in \mathbb{R}^n, \quad u(t) \in \USet, \quad t \in [t_0, \vartheta].
\end{equation}
Here, $t$ is the time;
$x(t)$ and $u(t)$ are respectively the current values of the state and control vectors;
the set $\USet \subset \mathbb{R}^r$ is compact, $r \in \mathbb{N}$;
$t_0$ and $\vartheta$ are the initial and terminal times.
We assume that the right-hand side of the differential equation in (\ref{system}) satisfies the following conditions:
\begin{description}
    \item[$(A.1)$]
        The function $A: [t_0, \vartheta] \rightarrow \mathbb{R}^{n \times n}$ is measurable and essentially bounded.

    \item[$(A.2)$]
        The function $f: [t_0, \vartheta] \times \USet \rightarrow \mathbb{R}^n$ is continuous.
\end{description}

Let $R_x > 0$ be fixed.
We suppose that, at the initial time $t_0$, an initial value $x_0 \in B(R_x) = \{y \in \mathbb{R}^n: \|y\| \leq R_x\}$ of the state vector of system (\ref{system}) is given.
As the set $\U(t_0, \vartheta)$ of admissible controls $u(\cdot)$ on the time interval $[t_0, \vartheta)$, we consider the set of all measurable functions $u: [t_0, \vartheta) \rightarrow \USet$.
By a motion of the system that corresponds to the initial value $x_0$ and a control $u(\cdot) \in \U(t_0, \vartheta)$, we mean a function $x(\cdot) \in \AC([t_0, \vartheta], \mathbb{R}^n)$ that satisfies the initial condition
\begin{equation} \label{initial_condition}
    x(t_0)
    = x_0
\end{equation}
and, together with $u(\cdot)$, satisfies the differential equation in (\ref{system}) for almost every $t \in [t_0, \vartheta]$.
Due to conditions $(A.1)$ and $(A.2)$, such a motion $x(\cdot)$ exists and is unique (see, e.g., \cite[Theorem~3.1]{Gomoyunov_2019_FCAA}), and we denote it by $x(\cdot \mid t_0, x_0, \vartheta, u(\cdot))$.

The goal of control is to minimize the cost functional
\begin{equation} \label{cost_functional}
    J(t_0, x_0, u(\cdot))
    = \sigma(x(\vartheta))
    + \int_{t_0}^{\vartheta} \chi(t, u(t)) \, \rd \tau,
    \quad u(\cdot) \in \U(t_0, \vartheta),
\end{equation}
where $x(\cdot) = x(\cdot \mid t_0, x_0, \vartheta, u(\cdot))$.
We assume that the conditions below hold:
\begin{description}
    \item[$(A.3)$]
        The function $\sigma: \mathbb{R}^n \rightarrow \mathbb{R}$ is continuous.

    \item[$(A.4)$]
        The function $\chi: [t_0, \vartheta] \times \USet \rightarrow \mathbb{R}$ is continuous.
\end{description}

The value of the optimal result in the control problem for system (\ref{system}) with initial condition (\ref{initial_condition}) and cost functional (\ref{cost_functional}) is defined by
\begin{equation} \label{Val}
    \Val(t_0, x_0)
    = \inf_{u(\cdot) \in \U(t_0, \vartheta)} J(t_0, x_0, u(\cdot)).
\end{equation}
A control $u^\circ(\cdot) \in \U(t_0, \vartheta)$ is called optimal in this problem, if the equality $J(t_0, x_0, u^\circ(\cdot)) = \Val(t_0, x_0)$ holds.

In order to find such an optimal control $u^\circ(\cdot)$, we propose to reduce the problem (\ref{system})--(\ref{cost_functional}) to some auxiliary optimal control problem in a dynamical system which motion is described by a first-order ordinary differential equation.
This reduction is based on some linear transformation, which is called the information image of a position of system (\ref{system}).

\subsection{Positions of the System}

According to \cite{Gomoyunov_2019_SIAM}, by a position of system (\ref{system}), we mean a pair $(t, w(\cdot))$ consisting of a time $t \in [t_0, \vartheta]$ and a function $w(\cdot) \in \AC([t_0, t], \mathbb{R}^n)$, $\|w(t_0)\| \leq R_x$, which is treated as a history of a motion of the system on the time interval $[t_0, t]$.
The set of all such positions is denoted by $G$.
Respectively, for every $x_0 \in B(R_x)$, the pair $(t_0, x_0) \in G$ is regarded as an initial position.

By analogy with Sect.~\ref{subsection_OCP_Statement_of_the_problem}, let us give a definition of motions of the system starting from an arbitrary position $(t_\ast, w_\ast(\cdot)) \in G$ and evolving on some time interval $[t_\ast, t^\ast]$, where $t^\ast \in [t_\ast, \vartheta]$.
Let us suppose that $t_\ast < t^\ast$.
Then, as the set $\U(t_\ast, t^\ast)$ of admissible controls $u(\cdot)$ on $[t_\ast, t^\ast)$, we consider the set of all measurable functions $u: [t_\ast, t^\ast) \rightarrow \USet$.
A motion of the system generated from the position $(t_\ast, w_\ast(\cdot))$ by a control $u(\cdot) \in \U(t_\ast, t^\ast)$ is defined as a function $x(\cdot) \in \AC([t_0, t^\ast], \mathbb{R}^n)$ that satisfies the equality
\begin{equation} \label{w_ast}
    x(t)
    = w_\ast(t),
    \quad t \in [t_0, t_\ast],
\end{equation}
and, together with $u(\cdot)$, satisfies the differential equation in (\ref{system}) for almost every $t \in [t_\ast, t^\ast]$.
By the scheme from \cite[Proposition~2]{Gomoyunov_2019_DGAA}, one can prove that, owing to conditions $(A.1)$ and $(A.2)$, such a motion $x(\cdot) = x(\cdot \mid t_\ast, w_\ast(\cdot), t^\ast, u(\cdot))$ exists and is unique.
In the degenerate case when $t_\ast = t^\ast$, the motion $x(\cdot)$ is completely determined by (\ref{w_ast}), and there is no
need in considering controls $u(\cdot)$ and determining the set $\U(t_\ast, t^\ast)$.
However, it is convenient to formally say that this motion $x(\cdot) = x(\cdot \mid t_\ast, w_\ast(\cdot), t^\ast, u(\cdot))$ is generated from $(t_\ast, w_\ast(\cdot))$ by $u(\cdot) \in \U(t_\ast, t^\ast)$.
Further, for the motion $x(\cdot)$ and a time $t \in [t_0, t^\ast]$, we denote the corresponding position of the system by $(t, x_t(\cdot))$, where the function $x_t: [t_0, t] \rightarrow \mathbb{R}^n$ is defined by
\begin{equation} \label{x_t}
    x_t(\tau)
    = x(\tau),
    \quad \tau \in [t_0, t].
\end{equation}
Let us note that the inclusion $(t, x_t(\cdot)) \in G$ is valid.

Following \cite{Gomoyunov_2019_FCAA_2} (see also \cite{Bourdin_2018}), let us consider the fundamental solution matrix of the differential equation in (\ref{system}), which is a continuous function
\begin{equation} \label{Omega}
    \Omega \ni (t, \tau) \mapsto F(t, \tau) \in \mathbb{R}^{n \times n},
    \quad \Omega
    = \{ (t, \tau) \in [t_0, \vartheta] \times [t_0, \vartheta]: t \geq \tau\},
\end{equation}
such that, for every fixed $\tau \in [t_0, \vartheta]$, the function $[\tau, \vartheta] \ni t \mapsto F(t, \tau) \in \mathbb{R}^{n \times n}$ is a unique continuous solution to the integral equation
\begin{equation*}
    F(t, \tau)
    = \frac{\Id}{\Gamma(\alpha)}
    + \frac{(t - \tau)^{1 - \alpha}}{\Gamma(\alpha)} \int_{\tau}^{t} \frac{A(\xi) F(\xi, \tau)}{(t - \xi)^{1 - \alpha} (\xi - \tau)^{1 - \alpha}} \, \rd \xi,
    \quad t \in [\tau, \vartheta].
\end{equation*}
Due to \cite[Theorem~5.2]{Gomoyunov_2019_FCAA_2}, for every motion $x(\cdot) = x(\cdot \mid t_\ast, w_\ast(\cdot), t^\ast, u(\cdot))$ of system (\ref{system}), where $(t_\ast, w_\ast(\cdot)) \in G$, $t^\ast \in [t_\ast, \vartheta]$, and $u(\cdot) \in \U(t_\ast, t^\ast)$, the representation formula below holds:
\begin{multline} \label{Cauchy_formula}
    x(t)
    = \Big( \Id + \int_{t_\ast}^{t} \frac{F(t, \tau) A(\tau)}{(t - \tau)^{1 - \alpha}} \, \rd \tau \Big) w_\ast(t_\ast) \\
    + \frac{1}{\Gamma(1 - \alpha)}\int_{t_\ast}^{t} \frac{F(t, \tau)}{(t - \tau)^{1 - \alpha}}
    \Big( \alpha \int_{t_0}^{t_\ast} \frac{w_\ast(\xi) - w_\ast(t_0)}{(\tau - \xi)^{1 + \alpha}} \, \rd \xi
    - \frac{w_\ast(t_\ast) - w_\ast(t_0)}{(\tau - t_\ast)^\alpha} \Big) \, \rd \tau \\
    + \int_{t_\ast}^{t} \frac{F(t, \tau) f(\tau, u(\tau))}{(t - \tau)^{1 - \alpha}} \, \rd \tau,
    \quad t \in [t_\ast, t^\ast].
\end{multline}

\subsection{Informational Image}

For a position $(t_\ast, w_\ast(\cdot)) \in G$, let us define the value $\Image(t_\ast, w_\ast(\cdot)) \in \mathbb{R}^n$, called the informational image of $(t_\ast, w_\ast(\cdot))$, as follows.
Let us consider the linear homogenous fractional differential equation corresponding to (\ref{system})
\begin{equation} \label{system_homogeneous}
    (^C D^\alpha y)(t) = A(t) y(t), \quad y(t) \in \mathbb{R}^n, \quad t \in [t_\ast, \vartheta],
\end{equation}
with the initial condition
\begin{equation} \label{initial_condition_homogenous}
    y(t)
    = w_\ast(t), \quad t \in [t_0, t_\ast].
\end{equation}
By analogy with the above, there exists a unique solution to the Cauchy problem (\ref{system_homogeneous}) and (\ref{initial_condition_homogenous}), which is the function $y(\cdot) = y(\cdot \mid t_\ast, w_\ast(\cdot), \vartheta) \in \AC([t_0, \vartheta], \mathbb{R}^n)$ that satisfies the equality in (\ref{initial_condition_homogenous}) and the differential equation in (\ref{system_homogeneous}) for almost every $t \in [t_\ast, \vartheta]$.
Then, we put
\begin{equation} \label{informational_image}
    \Image(t_\ast, w_\ast(\cdot))
    = y(\vartheta \mid t_\ast, w_\ast(\cdot), \vartheta).
\end{equation}
Let us note that, owing to representation formula (\ref{Cauchy_formula}), the informational image $\Image(t_\ast, w_\ast(\cdot))$ can be defined explicitly:
\begin{multline} \label{informational_image_Cauchy}
    \Image(t_\ast, w_\ast(\cdot))
    = \Big( \Id + \int_{t_\ast}^{\vartheta} \frac{F(\vartheta, \tau) A(\tau)}{(\vartheta - \tau)^{1 - \alpha}} \, \rd \tau \Big) w_\ast(t_\ast) \\
    + \frac{1}{\Gamma(1 - \alpha)}\int_{t_\ast}^{\vartheta} \frac{F(\vartheta, \tau)}{(\vartheta - \tau)^{1 - \alpha}}
    \Big( \alpha \int_{t_0}^{t_\ast} \frac{w_\ast(\xi) - w_\ast(t_0)}{(\tau - \xi)^{1 + \alpha}} \, \rd \xi
    - \frac{w_\ast(t_\ast) - w_\ast(t_0)}{(\tau - t_\ast)^\alpha} \Big) \, \rd \tau.
\end{multline}

Taking the last term in (\ref{Cauchy_formula}) into account, let us introduce the function
\begin{equation} \label{f_ast}
    f_\ast(t, u)
    = \frac{F(\vartheta, t) f(t, u)}{(\vartheta - t)^{1 - \alpha}},
    \quad t \in [t_0, \vartheta), \quad u \in \USet.
\end{equation}
Due to continuity of $F$ and $f$, the function $f_\ast: [t_0, \vartheta) \times \USet \rightarrow \mathbb{R}^n$ is continuous, and the following estimate is valid:
\begin{equation} \label{f_ast_bound}
    \|f_\ast(t, u)\|
    \leq \frac{M_F M_f}{(\vartheta - t)^{1 - \alpha}},
    \quad t \in [t_0, \vartheta), \quad u \in \USet,
\end{equation}
where we denote
\begin{equation}\label{M_F_M_f}
    M_F = \max_{t \in [t_0, \vartheta]} \|F(\vartheta, t)\|,
    \quad M_f = \max_{(t, u) \in [t_0, \vartheta] \times \USet} \|f(t, u)\|.
\end{equation}

The proposition below describes the dynamics of informational image (\ref{informational_image}) along motions of system (\ref{system}).
\begin{proposition} \label{Proposition_Informational_image}
    Let $x_0 \in B(R_x)$, $t^\ast \in [t_0, \vartheta]$, and $u(\cdot) \in \U(t_0, t^\ast)$.
    Let $x(\cdot) = x(\cdot \mid t_0, x_0, t^\ast, u(\cdot))$ be the corresponding motion of system $(\ref{system})$.
    Then, the equality below holds:
    \begin{equation} \label{Proposition_Informational_image_main}
        \Image(t, x_t(\cdot))
        = \Image(t_0, x_0)
        + \int_{t_0}^{t} f_\ast(\tau, u(\tau)) \, \rd \tau,
        \quad t \in [t_0, t^\ast],
    \end{equation}
    where $x_t(\cdot)$ is defined by $x(\cdot)$ according to $(\ref{x_t})$.
\end{proposition}

The proof of the proposition is given in Appendix.

\subsection{Auxiliary Optimal Control Problem}

Based on Proposition~\ref{Proposition_Informational_image}, we consider the auxiliary dynamical system which motion is described by the differential equation
\begin{equation} \label{system_z}
    \dot{z}(t) = f_\ast(t, p(t)),
    \quad z(t) \in \mathbb{R}^n, \quad p(t) \in \USet, \quad t \in [t_0, \vartheta),
\end{equation}
with the initial condition
\begin{equation} \label{initial_condition_z}
    z(t_0)
    = z_0
    = \Image(t_0, x_0).
\end{equation}
Here, $\dot{z}(t) = \rd z / \rd t$;
$z(t)$ and $p(t)$ are respectively the current values of the state and control vectors in the auxiliary system;
the set $\USet$ is the same as in original system (\ref{system});
the function $f_\ast$ is defined in (\ref{f_ast});
an initial value $z_0$ is determined by the informational image $\Image(t_0, x_0)$ (see (\ref{informational_image})) of an initial position $(t_0, x_0) \in G$ of system (\ref{system}).

By a position of auxiliary system (\ref{system_z}), we mean a pair $(t, z) \in [t_0, \vartheta] \times \mathbb{R}^n$.
Let $(t_\ast, z_\ast) \in [t_0, \vartheta] \times \mathbb{R}^n$, $t^\ast \in [t_\ast, \vartheta]$, and $p(\cdot) \in \U(t_\ast, t^\ast)$.
A motion $z(\cdot)$ of this system generated from $(t_\ast, z_\ast)$ by $p(\cdot)$ is an absolutely continuous function $z: [t_\ast, t^\ast] \rightarrow \mathbb{R}^n$ that satisfies the equality $z(t_\ast) = z_\ast$ and, together with $p(\cdot)$, satisfies the differential equation in (\ref{system_z}) for almost every $t \in [t_\ast, t^\ast]$.
In view of the described above properties of the function $f_\ast$, such a motion exists and is unique, and we denote it by $z(\cdot \mid t_\ast, z_\ast, t^\ast, p(\cdot))$.

As a direct consequence of Proposition~\ref{Proposition_Informational_image}, we derive the following lemma, giving a connection between motions of the original and auxiliary systems.
\begin{lemma} \label{Lemma_Connection}
    Let $x_0 \in B(R_x)$, $t^\ast \in [t_0, \vartheta]$, and $u(\cdot) \in \U(t_0, t^\ast)$.
    Let $x(\cdot) = x(\cdot \mid t_0, x_0, t^\ast, u(\cdot))$ be the corresponding motion of system $(\ref{system})$ with initial condition $(\ref{initial_condition})$.
    Let $z(\cdot) = z(\cdot \mid t_0, z_0, t^\ast, p(\cdot))$ be the motion of auxiliary system $(\ref{system_z})$ with initial condition $(\ref{initial_condition_z})$ generated by the same control $p(t) = u(t)$, $t \in [t_0, t^\ast)$.
    Then, the equality below holds:
    \begin{equation*}
        \Image(t, x_t(\cdot))
        = z(t),
        \quad t \in [t_0, t^\ast].
    \end{equation*}
\end{lemma}

Further, taking into account that (see (\ref{informational_image}))
\begin{equation} \label{informational_image_terminal}
    \Image(\vartheta, w_\ast(\cdot))
    = w_\ast(\vartheta),
    \quad (\vartheta, w_\ast(\cdot)) \in G,
\end{equation}
we define the auxiliary cost functional to be minimized by $p(\cdot)$ as follows:
\begin{equation} \label{cost_functional_z}
    J_\ast(t_0, z_0, p(\cdot))
    = \sigma(z(\vartheta))
    + \int_{t_0}^{\vartheta} \chi(\tau, p(\tau)) \, \rd \tau,
    \quad p(\cdot) \in \U(t_0, \vartheta),
\end{equation}
where $z(\cdot) = z(\cdot \mid t_0, z_0, \vartheta, p(\cdot))$, and the functions $\sigma$ and $\chi$ are taken from (\ref{cost_functional}).

The value of the optimal result in the auxiliary control problem for system (\ref{system_z}) with initial condition (\ref{initial_condition_z}) and cost functional (\ref{cost_functional_z}) is given by
\begin{equation*}
    \Vala(t_0, z_0)
    = \inf_{p(\cdot) \in \U(t_0, \vartheta)} J_\ast(t_0, z_0, p(\cdot)),
\end{equation*}
and a control $p^\circ(\cdot) \in \U(t_0, \vartheta)$ is called optimal if $J_\ast(t_0, z_0, p^\circ(\cdot)) = \Vala(t_0, z_0)$.

\begin{theorem} \label{Theorem_OCP}
    For any $x_0 \in B(R_x)$, the following statements are valid:
    \begin{description}
        \item[$i)$]
            A control $u^\circ(\cdot) \in \U(t_0, \vartheta)$ is optimal in the problem $(\ref{system})$--$(\ref{cost_functional})$ if and only if it is optimal in the auxiliary problem $(\ref{system_z})$, $(\ref{initial_condition_z})$, and $(\ref{cost_functional_z})$.
        \item[$ii)$]
            The optimal results in the original and auxiliary problems coincide, i.e., $\Val(t_0, x_0) = \Vala(t_0, z_0)$.
    \end{description}
\end{theorem}

The proof of the theorem is given in Appendix.

Thus, the original optimal control problem in fractional-order system (\ref{system}) is reduced to the auxiliary optimal control problem in first-order system (\ref{system_z}).

\begin{remark} \label{Remark_Reduction}
    Let us assume that the function $\sigma$ from cost functional (\ref{cost_functional}) can be represented in the following form:
    \begin{equation} \label{Remark_Reduction_sigma}
        \sigma(x)
        = \mu( K (x - c)),
        \quad x \in \mathbb{R}^n.
    \end{equation}
    Here, $K$ is a $(n_\ast \times n)$-matrix, $n_\ast \in \mathbb{N}$, $n_\ast < n$;
    $c \in \mathbb{R}^n$;
    $\mu: \mathbb{R}^{n_\ast} \rightarrow \mathbb{R}$ is a continuous function.
    In particular, this is the case when $\sigma$ does not depend on some $n - n_\ast$ coordinates of $x$.
    Under this additional assumption, we can reduce the dimension of the state vector in the auxiliary system from $n$ to $n_\ast$ and simplify the auxiliary cost functional.
    Namely, we consider the system
    \begin{equation*}
        \dot{z}(t) = K f_\ast(t, p(t)),
        \quad z(t) \in \mathbb{R}^{n_\ast}, \quad p(t) \in \USet, \quad t \in [t_0, \vartheta),
    \end{equation*}
    with the initial condition $z(t_0) = z_0 = K (\Image(t_0, x_0) - c)$ and the cost functional
    \begin{equation*}
        J_\ast(t_0, z_0, p(\cdot))
        = \mu(z(\vartheta))
        + \int_{t_0}^{\vartheta} \chi(t, p(t)) \, \rd t,
        \quad p(\cdot) \in \U(t_0, \vartheta).
    \end{equation*}
    One can show that, for this auxiliary optimal control problem, the result similar to Theorem~\ref{Theorem_OCP} takes place.
\end{remark}

However, it should be noted that, in general, the right-hand side of the differential equation in (\ref{system_z}) is unbounded, and, therefore, the auxiliary problem does not satisfy the assumptions that are usually made in the optimal control theory.
This complicates the application of the methods and results developed within this theory to solving the original optimal control problem.
In order to overcome this difficulty, in Sect.~\ref{subsection_OCP_Auxiliary_Optimal_Control_Problem_with_Parameter}, we propose to introduce a small parameter $\eta \in (0, \vartheta - t_0)$ and consider auxiliary system (\ref{system_z}) only up to the shifted terminal time $\vartheta_\eta = \vartheta - \eta$.
But before doing this, let us give an example illustrating Theorem~\ref{Theorem_OCP} and Remark~\ref{Remark_Reduction}.

\begin{example}
    Following \cite[Sect.~4.2]{Bergounioux_Bourdin_2019}, let us consider the optimal control problem described by the system
    \begin{equation} \label{Example_1_system}
        \begin{array}{c}
            \left\{
                \begin{array}{l}
                    (^C D^\alpha x_1)(t) = x_2(t) + \cos(u(t)), \\[0.5em]
                    (^C D^\alpha x_2)(t) = \sin(u(t)),
                \end{array}
            \right. \\[1.5em]
            x(t) = (x_1(t), x_2(t)) \in \mathbb{R}^2,
            \quad u(t) \in [- \pi / 2, \pi / 2],
            \quad t \in [t_0, \vartheta],
        \end{array}
    \end{equation}
    with the initial condition $x(t_0) = x_0 = (0, 0)$ and the cost functional
    \begin{equation*}
        J(t_0, x_0, u(\cdot))
        = - x_1(\vartheta),
        \quad u(\cdot) \in \U(t_0, \vartheta).
    \end{equation*}

    According to \cite[Theorem~4.2]{Idczak_Kamocki_2011} (see also \cite{Bourdin_2018,Gomoyunov_2019_FCAA_2}), the fundamental solution matrix $F$ of the differential equation in (\ref{Example_1_system}) is given by
    \begin{equation} \label{Example_1_F}
        F(t, \tau)
        = \sum_{i = 0}^\infty \frac{(t - \tau)^{i \alpha} A^i}{\Gamma((i + 1) \alpha)}
        = \begin{pmatrix}
            \frac{1}{\Gamma(\alpha)} & \frac{(t - \tau)^\alpha}{\Gamma(2 \alpha)} \\
            0 & \frac{1}{\Gamma(\alpha)}
        \end{pmatrix},
        \quad A
        = \begin{pmatrix}
            0 & 1 \\
            0 & 0
        \end{pmatrix},
    \end{equation}
    where $(t, \tau) \in \Omega$ (see (\ref{Omega})).
    Further, due to (\ref{informational_image_Cauchy}), we have
    \begin{equation} \label{Example_1_I}
        \Image(t_0, x_0) = (0, 0).
    \end{equation}
    Thus, denoting
    \begin{equation} \label{Example_1_b}
        b_1(t)
        = \frac{1}{\Gamma(\alpha) (\vartheta - t)^{1 - \alpha}},
        \quad
        b_2(t)
        = \frac{1}{\Gamma(2 \alpha) (\vartheta - t)^{1 - 2 \alpha}},
        \quad t \in [t_0, \vartheta),
    \end{equation}
    we come to the auxiliary optimal control problem for the system
    \begin{equation*}
        \begin{array}{c}
            \displaystyle
            \dot{z}(t)
            = b_1(t) \cos(p(t)) + b_2(t) \sin(p(t)), \\[0.5em]
            \displaystyle
            z(t) \in \mathbb{R}, \quad p(t) \in [- \pi / 2, \pi / 2], \quad t \in [t_0, \vartheta),
        \end{array}
    \end{equation*}
    with the initial condition $z(t_0) = z_0 = 0$ and the cost functional
    \begin{equation*}
        J_\ast(t_0, z_0, p(\cdot))
        = - z(\vartheta),
        \quad p(\cdot) \in \U(t_0, \vartheta).
    \end{equation*}

    By direct calculations, we obtain
    \begin{multline*}
        \Vala(t_0, z_0)
        = - \sup_{p(\cdot) \in \U(t_0, \vartheta)} \int_{t_0}^{\vartheta}
        \big( b_1(t) \cos(p(t)) + b_2(t) \sin(p(t)) \big) \, \rd t \\
        = - \int_{t_0}^{\vartheta}
        \big( b_1(t) \cos(p^\circ(t)) + b_2(t) \sin(p^\circ(t)) \big) \, \rd t
    \end{multline*}
    for the control $p^\circ(t) = \arctan( \Gamma(\alpha) (\vartheta - t)^\alpha / \Gamma(2 \alpha))$, $t \in [t_0, \vartheta]$.
    Then, this control $p^\circ(\cdot)$ is optimal in the auxiliary problem.
    Hence, by Theorem~\ref{Theorem_OCP} and Remark~\ref{Remark_Reduction}, we conclude that $p^\circ(\cdot)$ is an optimal control in the original problem, too.
    Let us note that this result agrees with \cite[Sect.~4.2]{Bergounioux_Bourdin_2019}.
\end{example}

\subsection{Auxiliary Optimal Control Problem with Parameter}
\label{subsection_OCP_Auxiliary_Optimal_Control_Problem_with_Parameter}

Let us fix $\eta \in (0, \vartheta - t_0)$, put $\vartheta_\eta = \vartheta - \eta \in (t_0, \vartheta)$, and consider the auxiliary optimal control problem for the system
\begin{equation} \label{system_z_eta}
    \dot{z}(t) = f_\ast(t, p_\eta(t)),
    \quad z(t) \in \mathbb{R}^n, \quad p_\eta(t) \in \USet, \quad t \in [t_0, \vartheta_\eta],
\end{equation}
with the initial condition
\begin{equation} \label{initial_condition_z_eta}
    z(t_0)
    = z_0
    = \Image(t_0, x_0)
\end{equation}
and the cost functional
\begin{equation} \label{cost_functional_z_eta}
    J_\eta(t_0, z_0, p_\eta(\cdot))
    = \sigma(z(\vartheta_\eta))
    + \int_{t_0}^{\vartheta_\eta} \chi (t, p_\eta(t)) \, \rd t,
    \quad p_\eta(\cdot) \in \U(t_0, \vartheta_\eta).
\end{equation}
The value of the optimal result in this problem is defined by
\begin{equation} \label{Valaeta}
    \Valaeta(t_0, z_0)
    = \inf_{p_\eta(\cdot) \in \U(t_0, \vartheta_\eta)} J_\eta(t_0, z_0, p_\eta(\cdot)),
\end{equation}
and a control $p^\circ_\eta(\cdot) \in \U(t_0, \vartheta_\eta)$ is called optimal if $J_\eta(t_0, z_0, p_\eta^\circ(\cdot)) = \Valaeta(t_0, z_0)$.
Thus, the only difference from the auxiliary optimal control problem (\ref{system_z}), (\ref{initial_condition_z}), and (\ref{cost_functional_z}) is that now the time $\vartheta_\eta$ is treated as the terminal one.

Let us note that, since the function (see (\ref{f_ast}))
\begin{equation} \label{f_ast_continuity}
    [t_0, \vartheta_\eta] \times \USet \ni (t, u) \mapsto f_\ast(t, u) \in \mathbb{R}^n
\end{equation}
is continuous (and, therefore, bounded), we obtain that, compared to system (\ref{system_z}), system (\ref{system_z_eta}) meets the typical assumptions from the optimal control theory.
On the other hand, in contrast to Theorem~\ref{Theorem_OCP}, due to the presence of the parameter $\eta$, one can not expect that an optimal control $p^\circ_\eta(\cdot)$ in the auxiliary problem (\ref{system_z_eta})--(\ref{cost_functional_z_eta}) determines some optimal control $u^\circ(\cdot)$ in the original problem (\ref{system})--(\ref{cost_functional}).

In this connection, let us consider a notion of $\varepsilon$-optimal controls.
Namely, for a number $\varepsilon > 0$, a control $u^\ast(\cdot) \in \U(t_0, \vartheta)$ is called $\varepsilon$-optimal in the original problem if $J(t_0, x_0, u^\ast(\cdot)) \leq \Val(t_0, x_0) + \varepsilon$.
Respectively, $p^\ast_\eta(\cdot) \in \U(t_0, \vartheta_\eta)$ is $\varepsilon$-optimal in the auxiliary problem if $J_\eta(t_0, z_0, p^\ast_\eta(\cdot)) \leq \Valaeta(t_0, z_0) + \varepsilon$.

\begin{theorem} \label{Theorem_OCP_varepsilon}
    For any $\varepsilon > 0$, there exist $\eta^\ast = \eta^\ast(\varepsilon) \in (0, \vartheta - t_0)$ and $\varepsilon^\ast = \varepsilon^\ast(\varepsilon) > 0$ such that, for any $\eta \in (0, \eta^\ast]$ and any $x_0 \in B(R_x)$, the following statements hold:
    \begin{description}
        \item[$i)$]
            If a control $p^\ast_\eta(\cdot)$ is $\varepsilon^\ast$-optimal in the auxiliary problem $(\ref{system_z_eta})$--$(\ref{cost_functional_z_eta})$ corresponding to the chosen $\eta$, then the control
            \begin{equation} \label{Theorem_OCP_varepsilon_main}
                u^\ast(t)
                = \begin{cases}
                    p^\ast_\eta(t), & \mbox{if } t \in [t_0, \vartheta_\eta), \\
                    \bar{u}, & \mbox{if } t \in (\vartheta_\eta, \vartheta),
                  \end{cases}
                \quad \bar{u} \in \USet,
            \end{equation}
            is $\varepsilon$-optimal in the problem $(\ref{system})$--$(\ref{cost_functional})$.
        \item[$ii)$]
            The optimal results in the original and auxiliary problems satisfy the inequality $|\Val(t_0, x_0) - \Valaeta(t_0, z_0)| \leq \varepsilon$.
    \end{description}
\end{theorem}

The proof of the theorem is given in Appendix.

Theorem~\ref{Theorem_OCP_varepsilon} allows us to apply, via the auxiliary problem (\ref{system_z_eta})--(\ref{cost_functional_z_eta}), the methods and results of the optimal control theory for first-order systems to finding $\varepsilon$-optimal controls in the problem (\ref{system})--(\ref{cost_functional}).
Let us note that a remark similar to Remark~\ref{Remark_Reduction} can also be made in relation to Theorem~\ref{Theorem_OCP_varepsilon}.
Let us consider an example.

\begin{example}
    Let the optimal control problem be described by the system
    \begin{equation} \label{Example_2_system}
        \begin{array}{c}
            \left\{
                \begin{array}{l}
                    (^C D^\alpha x_1)(t) = x_2(t), \\[0.5em]
                    (^C D^\alpha x_2)(t) = u(t),
                \end{array}
            \right. \\[1.5em]
            x(t) = (x_1(t), x_2(t)) \in \mathbb{R}^2,
            \quad u(t) \in [- 1, 1],
            \quad t \in [t_0, \vartheta],
        \end{array}
    \end{equation}
    with the initial condition $x(t_0) = x_0 = (0, 0)$ and the cost functional
    \begin{equation*}
        J(t_0, x_0, u(\cdot))
        = (x_1(\vartheta) - c_1)^2 + \int_{t_0}^{\vartheta} u^2(t) \, \rd t,
        \quad u(\cdot) \in \U(t_0, \vartheta),
    \end{equation*}
    where $c_1 \in \mathbb{R}$ is a given number.

    For every $\eta \in (0, \vartheta - t_0)$, taking (\ref{Example_1_F})--(\ref{Example_1_b}) into account, we come to the auxiliary optimal control problem for the system
    \begin{equation} \label{Example_2_system_auxiliary}
        \dot{z}(t)
        = b_2(t) p_\eta(t),
        \quad z(t) \in \mathbb{R}, \quad p_\eta(t) \in [- 1, 1], \quad t \in [t_0, \vartheta_\eta],
    \end{equation}
    with the initial condition $z(t_0) = z_0 = - c_1$ and the cost functional
    \begin{equation*}
        J_\eta(t_0, z_0, p_\eta(\cdot))
        = z^2(\vartheta_\eta) + \int_{t_0}^{\vartheta_\eta} p_\eta^2(t) \, \rd t,
        \quad p_\eta(\cdot) \in \U(t_0, \vartheta_\eta).
    \end{equation*}
    Applying the Pontryagin maximum principle to this auxiliary problem, we obtain (see, e.g., \cite[Ch.~3, Theorem~14]{Lee_Markus_1967} and also \cite[Statement~1]{Shaburov_2017}) that the unique optimal control is given by $p^\circ_\eta(t) = b_2(t) \lambda_\eta / S(b_2(t) |\lambda_\eta|)$, $t \in [t_0, \vartheta_\eta)$, where $S(\zeta) = 2$ for $\zeta \in [0, 2]$ and $S(\zeta) = \zeta$ for $\zeta > 2$, and $\lambda_\eta \in \mathbb{R}$ is the unique solution to the equation
    \begin{equation*}
        \lambda \int_{t_0}^{\vartheta_\eta} \frac{b_2^{\, 2}(t)}{S(b_2(t) |\lambda|)} \, \rd t + \frac{\lambda}{2}
        = c_1.
    \end{equation*}

    Consequently, by Theorem~\ref{Theorem_OCP_varepsilon} and Remark~\ref{Remark_Reduction}, for every $\varepsilon > 0$, on the basis of the control $p^\circ_\eta(\cdot)$ for a sufficiently small $\eta \in (0, \vartheta - t_0)$, we can determine an $\varepsilon$-optimal in the original problem control $u^\ast(\cdot)$ according to (\ref{Theorem_OCP_varepsilon_main}).
\end{example}

\section{Optimal Feedback Controls}
\label{section_OFC}

In the previous section, we study the question of finding optimal and $\varepsilon$-optimal (open-loop) controls $u(\cdot) \in \U(t_0, \vartheta)$ in the optimal control problem (\ref{system})--(\ref{cost_functional}).
However, it is often more convenient to use feedback (closed-loop) controls, since, in many cases, they are easier to construct, and, moreover, they do not depend on a particular choice of an initial value $x_0 \in B(R_x)$.
In this section, we develop the proposed above reduction of the original problem to the auxiliary problem (\ref{system_z_eta})--(\ref{cost_functional_z_eta}) for obtaining optimal feedback controls.

\subsection{Positional Control Strategies}

We consider a formalization of feedback controls within the framework of positional control strategies \cite{Krasovskii_Subbotin_1988,Krasovskii_Krasovskii_1995} (see also \cite{Gomoyunov_2019_Trudy_Eng,Gomoyunov_2019_SIAM}).
By a (positional) control strategy, we mean an arbitrary function
\begin{equation*}
    G \times (0, \vartheta - t_0) \ni (t, w(\cdot), \eta)
    \mapsto U(t, w(\cdot), \eta) \in \USet,
\end{equation*}
where $\eta$ plays a role of some accuracy parameter.

Let us fix $\eta \in (0, \vartheta - t_0)$ and a partition $\Delta$ of the time interval $[t_0, \vartheta]$:
\begin{equation} \label{Delta}
    \Delta = \{\tau_j\}_{j = 1}^{k + 1},
    \quad \tau_1 = t_0, \quad \tau_{j+1} > \tau_j, \quad j \in \overline{1, k}, \quad \tau_{k+1} = \vartheta, \quad k \in \mathbb{N}.
\end{equation}
The triple $\{U, \eta, \Delta\}$ is called a control law.
This control law forms in system (\ref{system}) a piecewise constant control $u(\cdot) \in \U(t_0, \vartheta)$ by the following feedback rule:
\begin{equation} \label{Control_law}
    u(t) = U(\tau_j, x_{\tau_j}(\cdot), \eta),
    \quad t \in [\tau_j, \tau_{j + 1}), \quad j \in \overline{1, k},
\end{equation}
where, as usual, we denote $x_{\tau_j}(t) = x(t)$, $t \in [t_0, \tau_j]$.
Let us note that, for every initial value $x_0 \in B(R_x)$, the control law $\{U, \eta, \Delta\}$ determines the control $u(\cdot) = u(\cdot \mid t_0, x_0, \vartheta, U, \eta, \Delta)$ and the corresponding motion $x(\cdot)$ uniquely.
Let us note also that, in accordance with (\ref{Val}), we have
\begin{equation*}
    J\big( t_0, x_0, u(\cdot \mid t_0, x_0, \vartheta, U, \eta, \Delta) \big)
    \geq \Val(t_0, x_0).
\end{equation*}

Taking this into account, we call a control strategy $U^\circ$ optimal (uniformly with respect to initial values $x_0 \in B(R_x)$), if the following statement holds.
For any $\varepsilon > 0$, there exist
\begin{equation*}
    \eta^\circ = \eta^\circ(\varepsilon) \in (0, \vartheta - t_0),
    \quad (0, \eta^\circ] \ni \eta \mapsto \delta^\circ(\eta) = \delta^\circ(\varepsilon, \eta) \in (0, \infty)
\end{equation*}
such that, for any $\eta \in (0, \eta^\circ]$, any partition $\Delta$ (\ref{Delta}) with the diameter $\diam(\Delta) = \max_{j \in \overline{1, k}} (\tau_{j+1} - \tau_j) \leq \delta^\circ(\eta)$, and any $x_0 \in B(R_x)$, the inequality
\begin{equation} \label{optimal_control_strategy}
    J\big( t_0, x_0, u(\cdot \mid t_0, x_0, \vartheta, U^\circ, \eta, \Delta) \big)
    \leq \Val(t_0, x_0) + \varepsilon
\end{equation}
is valid, i.e., the control $u(\cdot \mid t_0, x_0, \vartheta, U^\circ, \eta, \Delta)$ is $\varepsilon$-optimal.

Following the ideas from Sect.~\ref{subsection_OCP_Auxiliary_Optimal_Control_Problem_with_Parameter}, let us construct such an optimal control strategy $U^\circ$ on the basis of optimal control strategies $P^\circ_\eta$ in the auxiliary problems (\ref{system_z_eta})--(\ref{cost_functional_z_eta}) for $\eta \in (0, \vartheta - t_0)$.

\subsection{Positional Control Strategies in the Auxiliary Problem}

Let us fix $\eta \in (0, \vartheta - t_0)$, consider the auxiliary optimal control problem (\ref{system_z_eta})--(\ref{cost_functional_z_eta}), and define the set $G_\eta = [t_0, \vartheta_\eta] \times \mathbb{R}^n$ of all positions of system (\ref{system_z_eta}).
A (positional) control strategy is a function
\begin{equation*}
    G_\eta \times (0, \vartheta - t_0) \ni (t, z, \varkappa) \mapsto P_\eta(t, z, \varkappa) \in \USet.
\end{equation*}
where $\varkappa$ is treated as an accuracy parameter.

Let $\varkappa \in (0, \vartheta - t_0)$, and let $\Delta_\eta$ be a partition of the time interval $[t_0, \vartheta_\eta]$:
\begin{equation} \label{Delta_eta}
    \Delta_\eta = \{\tau_j\}_{j = 1}^{k + 1},
    \quad \tau_1 = t_0, \quad \tau_{j+1} > \tau_j, \quad j \in \overline{1, k}, \quad \tau_{k + 1} = \vartheta_\eta, \quad k \in \mathbb{N}.
\end{equation}
The control law $\{P_\eta, \varkappa, \Delta_\eta\}$ forms in the auxiliary system a piecewise constant control $p_\eta(\cdot) \in \U(t_0, \vartheta_\eta)$ by the following feedback rule:
\begin{equation} \label{Control_law_eta}
    p_\eta(t)
    = P_\eta(\tau_j, z(\tau_j), \varkappa),
    \quad t \in [\tau_j, \tau_{j + 1}), \quad j \in \overline{1, k}.
\end{equation}
For every $z_0 \in \mathbb{R}^n$, the control law $\{P_\eta, \varkappa, \Delta_\eta\}$ determines the control $p_\eta(\cdot) = p_\eta(\cdot \mid t_0, z_0, \vartheta_\eta, P_\eta, \varkappa, \Delta_\eta)$ and the corresponding motion $z(\cdot)$ uniquely.

In accordance with condition $(A.1)$, taking $M_F$ from (\ref{M_F_M_f}), let us denote
\begin{equation} \label{R_z}
    M_A
    = \esssup{t \in [t_0, \vartheta]} \|A(t)\|,
    \quad R_z
    = \big(1 + M_F M_A (\vartheta - t_0)^\alpha / \alpha \big) R_x.
\end{equation}
Then, due to (\ref{informational_image_Cauchy}), for every initial value $x_0 \in B(R_x)$ from (\ref{initial_condition}), we have
\begin{equation} \label{I_x_0_R_z}
    \|\Image(t_0, x_0)\|
    \leq \Big\| \Big( \Id + \int_{t_0}^{\vartheta} \frac{F(\vartheta, \tau) A(\tau)}{(\vartheta - \tau)^{1 - \alpha}} \, \rd \tau \Big) x_0 \Big\|
    \leq R_z,
\end{equation}
Therefore, in view of (\ref{initial_condition_z_eta}), in the auxiliary problem, we can restrict ourselves to initial values $z_0 \in B(R_z)$.

Thus, we call a control strategy $P_\eta^\circ$ optimal (uniformly with respect to $z_0 \in B(R_z)$), if the following statement holds.
For any $\varepsilon > 0$, there exist
\begin{equation} \label{kappa_delta_eta_a}
    \varkappa_\eta^\circ = \varkappa_\eta^\circ(\varepsilon) \in (0, \vartheta - t_0),
    \quad (0, \varkappa_\eta^\circ] \ni \varkappa \mapsto \omega_\eta^\circ(\varkappa) = \omega_\eta^\circ(\varepsilon, \varkappa) \in (0, \infty)
\end{equation}
such that, for any $\varkappa \in (0, \varkappa_\eta^\circ]$, any partition $\Delta_\eta$ (\ref{Delta_eta}) with the diameter $\diam(\Delta_\eta) \leq \omega_\eta^\circ(\varkappa)$, and any $z_0 \in B(R_z)$, the inequality below is valid:
\begin{equation} \label{optimal_control_strategy_a}
    J_\eta\big( t_0, z_0, p_\eta(\cdot \mid t_0, z_0, \vartheta_\eta, P_\eta^\circ, \varkappa, \Delta_\eta) \big)
    \leq \Valaeta(t_0, z_0) + \varepsilon,
\end{equation}
i.e., $p_\eta(\cdot \mid t_0, z_0, \vartheta_\eta, P_\eta^\circ, \varkappa, \Delta_\eta)$ is an $\varepsilon$-optimal control in the auxiliary problem.

Let us note that, due to continuity of $f_\ast$ (see (\ref{f_ast_continuity})) and conditions $(A.3)$ and $(A.4)$, according to, e.g., \cite[Theorems~9.2 and~22.1]{Krasovskii_Krasovskii_1995} (see also \cite[Sect.~6]{Lukoyanov_Gomoyunov_2019_DGAA} and the references therein), such an optimal control strategy $P^\circ_\eta$ exists.
Moreover, it can be constructed, for example, by the method of extremal shift to accompanying points, which is shortly described below.

In the auxiliary problem, we consider the value function (see (\ref{cost_functional_z_eta}) and (\ref{Valaeta}))
\begin{equation*}
    \Valaeta(t, z)
    = \inf_{p_\eta(\cdot) \in \U(t, \vartheta_\eta)}
    \Big( \sigma(z(\vartheta_\eta)) + \int_{t}^{\vartheta_\eta} \chi(\tau, p_\eta(\tau))  \, \rd \tau \Big),
    \quad (t, z) \in G_\eta,
\end{equation*}
where $z(\cdot) = z(\cdot \mid t, z, \vartheta_\eta, p_\eta(\cdot))$.
Let $\varkappa \in (0, \vartheta - t_0)$ be fixed.
For every position $(t, z) \in G_\eta$, relying on the value function, we choose the accompanying point
\begin{equation} \label{OFC_extremal_shift_1}
    (z^\circ, z_{n + 1}^\circ)
    \in \argmin{(\bar{z}, \bar{z}_{n + 1})} (\Valaeta(t, \bar{z}) + \bar{z}_{n + 1}),
\end{equation}
where the minimum is calculated over the pairs $(\bar{z}, \bar{z}_{n + 1}) \in \mathbb{R}^n \times \mathbb{R}$ such that
\begin{equation} \label{OFC_extremal_shift_2}
    \|z - \bar{z}\|^2 + \bar{z}_{n + 1}^2
    \leq r^2(t, \varkappa),
    \quad r^2(t, \varkappa)
    = \varkappa + (t - t_0) \varkappa,
\end{equation}
and, after that, we determine
\begin{equation} \label{OFC_extremal_shift_3}
    P^\circ_\eta(t, z, \varkappa)
    \in \argmin{p \in \USet} \big( \langle z - z^\circ, f_\ast(t, p) \rangle - z_{n + 1}^\circ \chi(t, p) \big).
\end{equation}

Now, we define a control strategy $U^\ast$ in the original problem as follows.
For every $\eta \in (0, \vartheta - t_0)$, let us consider an optimal control strategy $P^\circ_\eta$ in the auxiliary problem and determine the corresponding number $\varkappa^\circ_\eta(\eta)$ according to (\ref{kappa_delta_eta_a}).
Then, for every $(t, w(\cdot)) \in G$, if $t < \vartheta_\eta$, we put
\begin{equation} \label{strategy_U_ast}
    U^\ast(t, w(\cdot), \eta)
    = P^\circ_\eta \big(t, \Image(t, w(\cdot)), \varkappa^\circ_\eta (\eta) \big),
\end{equation}
where $\Image(t, w(\cdot))$ is the informational image (see (\ref{informational_image}) or (\ref{informational_image_Cauchy})).
If $t \in [\vartheta_\eta, \vartheta]$, we formally define $U^\ast(t, w(\cdot), \eta) = \bar{u}$ for some fixed $\bar{u} \in \USet$.

\begin{theorem} \label{Theorem_OFC}
    The control strategy $U^\ast$ defined by $(\ref{strategy_U_ast})$ is optimal in the problem $(\ref{system})$--$(\ref{cost_functional})$.
\end{theorem}

The proof of the theorem is given in Appendix.

Thus, in order to construct an optimal control strategy $U^\circ$ in the original problem in fractional-order system (\ref{system}), it is sufficient to find for every $\eta \in (0, \vartheta - t_0)$ an optimal control strategy $P^\circ_\eta$ in the auxiliary problem in first-order system (\ref{system_z_eta}).

\begin{remark} \label{Remark_Reduction_FC}
    As in Remark~\ref{Remark_Reduction}, let us suppose that the function $\sigma$ from cost functional (\ref{cost_functional}) can be represented as in (\ref{Remark_Reduction_sigma}).
    Then, for every $\eta \in (0, \vartheta - t_0)$, the auxiliary optimal control problem is described by the system
    \begin{equation*}
        \dot{z}(t) = K f_\ast(t, p_\eta(t)),
        \quad z(t) \in \mathbb{R}^{n_\ast}, \quad p_\eta(t) \in \USet, \quad t \in [t_0, \vartheta_\eta],
    \end{equation*}
    and the cost functional
    \begin{equation*}
        J_\eta(t_0, z_0, p_\eta(\cdot))
        = \mu(z(\vartheta))
        + \int_{t_0}^{\vartheta_\eta} \chi(t, p_\eta(t)) \, \rd t,
        \quad p_\eta(\cdot) \in \U(t_0, \vartheta_\eta).
    \end{equation*}
    Let $P^\circ_\eta$ be an optimal control strategy in this auxiliary problem.
    In accordance with (\ref{strategy_U_ast}), we put
    \begin{equation} \label{strategy_U_ast_Remark_Reduction}
        U^\ast (t, w(\cdot), \eta)
        = P^\circ_\eta \big(t, K(\Image(t, w(\cdot)) - c), \varkappa^\circ_\eta (\eta) \big)
    \end{equation}
    for $(t, w(\cdot)) \in G$ such that $t < \vartheta_\eta$.
    By analogy with Theorem~\ref{Theorem_OFC}, one can prove that such a control strategy $U^\ast$ is optimal in the original problem.
\end{remark}

The following two examples illustrate Theorem~\ref{Theorem_OFC} and Remark~\ref{Remark_Reduction_FC}.

\begin{example}
    Let us consider the optimal control problem for system (\ref{Example_2_system}) and the cost functional
    \begin{equation*}
        J(t_0, x_0, u(\cdot))
        = (x_1(\vartheta) - c_1)^2,
        \quad u(\cdot) \in \U(t_0, \vartheta),
    \end{equation*}
    where $c_1 \in \mathbb{R}$ is a given number.
    Then, for every $\eta \in (0, \vartheta - t_0)$, we come to the auxiliary optimal control problem for system (\ref{Example_2_system_auxiliary}) and the cost functional
    \begin{equation*}
        J(t_0, z_0, p_\eta(\cdot))
        = z^2(\vartheta),
        \quad p_\eta(\cdot) \in \U(t_0, \vartheta_\eta).
    \end{equation*}
    Let us associate this auxiliary problem with the Cauchy problem for the corresponding Hamilton--Jacobi--Bellman equation
    \begin{equation*}
        \frac{\partial \varphi(t, z)}{\partial t} - b_2(t) \Big| \frac{\partial \varphi(t, z)}{\partial z} \Big|
        = 0,
        \quad (t, z) \in (t_0, \vartheta_\eta) \times \mathbb{R},
    \end{equation*}
    with the right-end boundary condition $\varphi(\vartheta_\eta, z) = z^2$, $z \in \mathbb{R}$.
    Let us define
    \begin{equation*}
        \varphi_\eta(t, z)
        = \begin{cases}
            (|z| - \psi_\eta(t))^2, & \mbox{if } |z| > \psi_\eta(t), \\
            0, & \mbox{if } |z| \leq \psi_\eta(t),
          \end{cases}
          \quad (t, z) \in G_\eta = [t_0, \vartheta_\eta] \times \mathbb{R},
    \end{equation*}
    where, for the function $b_2$ from (\ref{Example_1_b}), we denote
    \begin{equation*}
        \psi_\eta(t)
        = \int_{t_0}^{\vartheta_\eta} b_2(t) \, \rd t
        = \frac{(\vartheta - t)^{2 \alpha} - \eta^{2 \alpha}}{\Gamma(2 \alpha + 1)},
        \quad t \in [t_0, \vartheta_\eta].
    \end{equation*}
    One can verify that $\varphi_\eta$ is a continuously differentiable solution to the considered Cauchy problem.
    Thus, according to, e.g., \cite[Theorem~4.1.1]{Krasovskii_Subbotin_1988}, we have $\Valaeta(t, z) = \varphi_\eta(t, z)$, $(t, z) \in G_\eta$, and the optimal control strategy in the auxiliary problem is given by
    \begin{equation*}
        P^\circ_\eta(t, z)
        \in \argmin{p \in [-1, 1]} \Big(p \frac{\partial \varphi_\eta(t, z)}{\partial z}\Big)
        = \begin{cases}
            \{1\}, & \mbox{if } z < - \psi_\eta(t), \\
            [-1, 1], & \mbox{if } |z| \leq \psi_\eta(t), \\
            \{- 1\}, & \mbox{if } z > \psi_\eta(t),
          \end{cases}
    \end{equation*}
    where $(t, z) \in [t_0, \vartheta_\eta) \times \mathbb{R}$.

    Hence, by Theorem~\ref{Theorem_OFC} and Remark~\ref{Remark_Reduction_FC}, on the basis of the found optimal control strategies $P^\circ_\eta$, $\eta \in (0, \vartheta - t_0)$, we can construct an optimal control strategy in the original problem as follows:
    \begin{equation*}
        U^\circ (t, w(\cdot), \eta)
        = \begin{cases}
            P^\circ_\eta \big(t, \Image_1(t, w(\cdot)) - c_1\big), & \mbox{if } t < \vartheta_\eta, \\
            \bar{u}, & \mbox{if } t \geq \vartheta_\eta,
        \end{cases}
    \end{equation*}
    where $(t, w(\cdot)) \in G$ and $\bar{u} \in \USet$, and $\Image_1(t, w(\cdot))$ is the first coordinate of the informational image $\Image(t, w(\cdot))$ defined by (\ref{informational_image}) or (\ref{informational_image_Cauchy}).
\end{example}

\begin{example}
    Let us consider the original optimal control problem (\ref{system})--(\ref{cost_functional}) in the case when $\sigma(x) = \|K (x - c)\|$, $x \in \mathbb{R}^n$ (see Remark~\ref{Remark_Reduction_FC}).
    Then, in the auxiliary optimal control problem (\ref{system_z_eta})--(\ref{cost_functional_z_eta}), the cost functional takes the form
    \begin{equation*}
        J_\eta(t_0, z_0, p_\eta(\cdot))
        = \|z(\vartheta_\eta)\|
        + \int_{t_0}^{\vartheta_\eta} \chi (t, p_\eta(t)) \, \rd t,
        \quad p_\eta(\cdot) \in \U(t_0, \vartheta_\eta).
    \end{equation*}
    Hence, according to, e.g., \cite[\S~23]{Krasovskii_Krasovskii_1995} (see also \cite[Sect.~7.1]{Lukoyanov_Gomoyunov_2019_DGAA} and the references therein), the value function in the auxiliary problem is given by
    \begin{equation*}
        \Valaeta (t, z)
        = \max_{l \in B(1)} \big( \langle l, z \rangle + \nu_\eta(t, l) \big),
        \quad (t, z) \in G_\eta = [t_0, \vartheta_\eta] \times \mathbb{R}^{n_\ast},
    \end{equation*}
    where $B(1) = \{l \in \mathbb{R}^{n_\ast}: \|l\| \leq 1\}$ and
    \begin{equation*}
        \nu_\eta(t, l)
        = \int_{t}^{\vartheta_\eta} \min_{p \in \USet}
        \big( \langle l, K f_\ast(\tau, p) \rangle + \chi(\tau, p) \big) \, \rd \tau,
        \quad t \in [t_0, \vartheta_\eta], \quad l \in B(1).
    \end{equation*}
    Moreover, applying the method of extremal shift to accompanying points (\ref{OFC_extremal_shift_1})--(\ref{OFC_extremal_shift_3}), one can construct an optimal control strategy as follows:
    \begin{equation*}
        P^\circ_\eta(t, z, \varkappa)
        \in \argmin{p \in \USet} \big( \langle l^\circ, K f_\ast(t, p) \rangle + \chi(t, p) \big),
        \quad (t, z) \in G_\eta, \quad \varkappa \in (0, \vartheta - t_0),
    \end{equation*}
    where
    \begin{equation*}
        l^\circ
        \in \argmax{l \in B(1)} \big( \langle l, z \rangle + \nu_\eta(t, l) - r(t, \varkappa) \sqrt{1 + \|l\|^2} \big).
    \end{equation*}
    Thus, based on these formulas and (\ref{strategy_U_ast_Remark_Reduction}), we can effectively calculate an optimal control strategy in the original problem.
\end{example}

\begin{acknowledgements}
    This work was supported by RSF, project no. 19-11-00105.
\end{acknowledgements}

\appendix
\section*{Appendix: Proofs}

{\it Proof of Proposition~\ref{Proposition_Informational_image}}
    Let us fix $x_0 \in B(R_x)$, $t^\ast \in [t_0, \vartheta]$, and $u(\cdot) \in \U(t_0, t^\ast)$, and consider the corresponding motion $x(\cdot) = x(\cdot \mid t_0, x_0, t^\ast, u(\cdot))$ of system (\ref{system}).
    Let us take $t_\ast \in [t_0, t^\ast]$ and prove the equality in (\ref{Proposition_Informational_image_main}) for $t = t_\ast$.
    According to (\ref{informational_image}), we have $\Image(t_0, x_0) = y_0(\vartheta)$, where $y_0(\cdot) = y(\cdot \mid t_0, x_0, \vartheta)$ is the solution to the following Cauchy problem:
    \begin{equation*}
        (^C D^\alpha y_0)(t)
        = A(t) y_0(t),
        \quad y_0(t) \in \mathbb{R}^n, \quad t \in [t_0, \vartheta];
        \quad y_0(t_0)
        = x_0.
    \end{equation*}
    Respectively, $\Image(t_\ast, x_{t_\ast}(\cdot)) = y_\ast(\vartheta)$, where $y_\ast(\cdot) = y(\cdot \mid t_\ast, x_{t_\ast}(\cdot), \vartheta)$ is the solution to
    \begin{equation*}
        (^C D^\alpha y_\ast)(t)
        = A(t) y_\ast(t),
        \quad y_\ast(t) \in \mathbb{R}^n, \quad t \in [t_\ast, \vartheta];
        \quad y_\ast(t)
        = x(t),
        \quad t \in [t_0, t_\ast].
    \end{equation*}
    Hence, for the difference $s(t) = y_\ast(t) - y_0(t)$, $t \in [t_0, \vartheta]$, we obtain $s(\cdot) \in \AC([t_0, \vartheta], \mathbb{R}^n)$, $s(t_0) = 0$, and, for almost every $t \in [t_0, \vartheta]$,
    \begin{equation*}
        (^C D^\alpha s)(t)
        = (^C D^\alpha y_\ast)(t) - (^C D^\alpha y_0)(t)
        = A(t) (y_\ast(t) - y_0(t)) + g(t)
        = A(t) s(t) + g(t),
    \end{equation*}
    where $g(t) = f(t, u(t))$ for $t \in [t_0, t_\ast)$ and $g(t) = 0$ for $t \in [t_\ast, \vartheta]$.
    Then, in view of representation formula (\ref{Cauchy_formula}) and notation (\ref{f_ast}), we get
    \begin{equation*}
        s(\vartheta)
        = \int_{t_0}^{\vartheta} \frac{F(\vartheta, \tau) g(\tau)}{(\vartheta - \tau)^{1 - \alpha}} \, \rd \tau
        = \int_{t_0}^{t_\ast} \frac{F(\vartheta, \tau) f(\tau, u(\tau))}{(\vartheta - \tau)^{1 - \alpha}} \, \rd \tau
        = \int_{t_0}^{t_\ast} f_\ast(\tau, u(\tau)) \, \rd \tau,
    \end{equation*}
    wherefrom, since $s(\vartheta) = \Image(t_\ast, x_{t_\ast}(\cdot)) - \Image(t_0, x_0)$, we derive the equality in (\ref{Proposition_Informational_image_main}) for $t = t_\ast$.
\qed

{\it Proof of Theorem~\ref{Theorem_OCP}}
    Let $x_0 \in B(R_x)$ be fixed, and let $z_0 = \Image(t_0, x_0)$.
    For every control $u(\cdot) \in \U(t_0, \vartheta)$, let us consider the corresponding motions $x(\cdot) = x(\cdot \mid t_0, x_0, \vartheta, u(\cdot))$ and $z(\cdot) = z(\cdot \mid t_0, z_0, \vartheta, u(\cdot))$ of original (\ref{system}) and auxiliary (\ref{system_z}) systems, respectively.
    Then, it follows from Lemma~\ref{Lemma_Connection} and (\ref{informational_image_terminal}) that $z(\vartheta) = \Image(\vartheta, x_\vartheta(\cdot)) = x(\vartheta)$, and, therefore,
    \begin{multline} \label{Theorem_OCP_proof_1}
        J(t_0, x_0, u(\cdot))
        = \sigma(x(\vartheta)) + \int_{t_0}^{\vartheta} \chi(t, u(t)) \, \rd t \\
        = \sigma(z(\vartheta)) + \int_{t_0}^{\vartheta} \chi(t, u(t)) \, \rd t
        = J_\ast(t_0, z_0, u(\cdot)).
    \end{multline}
    Since this equality holds for every $u(\cdot) \in \U(t_0, \vartheta)$, we get the statements of the theorem.
\qed

{\it Proof of Theorem~\ref{Theorem_OCP_varepsilon}}
    Let $\varepsilon > 0$ be fixed.
    Taking $M_F$ and $M_f$ from (\ref{M_F_M_f}) and $R_z$ from (\ref{R_z}), we define
    \begin{equation*}
        M_z
        = R_z + M_F M_f (\vartheta - t_0)^\alpha / \alpha.
    \end{equation*}
    Due to continuity of $\sigma$ (see condition $(A.3)$), there exists $\zeta > 0$ such that, for any $z_1$, $z_2 \in B(M_z)$, if $\|z_1 - z_2\| \leq \zeta$, then $|\sigma(z_1) - \sigma(z_2)| \leq \varepsilon / 6$.
    Let us choose $\eta_1 > 0$ such that $M_F M_f \eta_1^\alpha / \alpha \leq \zeta$.
    Further, denoting (see condition $(A.4)$)
    \begin{equation*}
        M_\chi
        = \max_{(t, u) \in [t_0, \vartheta] \times \USet} |\chi(t, u)|,
    \end{equation*}
    we take $\eta_2 > 0$ satisfying $M_\chi \eta_2 \leq \varepsilon / 6$ and $\eta_2 < \vartheta - t_0$ and put $\eta^\ast = \min\{\eta_1, \eta_2\}$.
    Let us show that the statements of the theorem are valid for the chosen $\eta^\ast$ and $\varepsilon^\ast = \varepsilon / 3$.

    We fix $\eta \in (0, \eta^\ast]$ and $x_0 \in B(R_x)$ and define $z_0 = \Image(t_0, x_0)$.
    Let $u(\cdot) \in \U(t_0, \vartheta)$, and let $z(\cdot) = z(\cdot \mid t_0, z_0, \vartheta, u(\cdot))$ be the motion of auxiliary system (\ref{system_z}).
    Due to (\ref{f_ast_bound}) and (\ref{system_z}), we obtain
    \begin{equation} \label{Theorem_OCP_varepsilon_proof_z_bound}
        \|z(t)\|
        \leq \|z_0\| + \int_{t_0}^{t} \|f_\ast(\tau, u(\tau))\| \, \rd \tau
        \leq R_z + M_F M_f (\vartheta - t_0)^\alpha / \alpha
        = M_z,
        \quad t \in [t_0, \vartheta],
    \end{equation}
    and
    \begin{equation} \label{Theorem_OCP_varepsilon_proof_z_Holder}
        \|z(\vartheta) - z(\vartheta_\eta)\|
        \leq \int_{\vartheta_\eta}^{\vartheta} \|f_\ast(\tau, u(\tau))\| \, \rd \tau
        \leq M_F M_f \eta^\alpha / \alpha
        \leq \zeta.
    \end{equation}
    Further, let $z_\eta(\cdot) = z(\cdot \mid t_0, z_0, \vartheta_\eta, p_\eta(\cdot))$ be the motion of auxiliary system (\ref{system_z_eta}) corresponding to the control $p_\eta(t) = u(t)$, $t \in [t_0, \vartheta_\eta)$.
    Then, we have $z_\eta(\vartheta_\eta) = z(\vartheta_\eta)$, and, owing to (\ref{Theorem_OCP_varepsilon_proof_z_bound}) and (\ref{Theorem_OCP_varepsilon_proof_z_Holder}), we derive
    \begin{multline*}
        |J_\ast(t_0, z_0, u(\cdot)) - J_\eta(t_0, z_0, p_\eta(\cdot))| \\
        \leq |\sigma(z(\vartheta)) - \sigma(z_\eta(\vartheta_\eta))|
        + \int_{\vartheta_\eta}^{\vartheta} |\chi(\tau, u(\tau))| \, \rd \tau
        \leq \varepsilon / 6 + M_\chi \eta
        \leq \varepsilon / 3.
    \end{multline*}
    Taking into account that $J(t_0, x_0, u(\cdot)) = J_\ast(t_0, z_0, u(\cdot))$ according to (\ref{Theorem_OCP_proof_1}), we conclude
    \begin{equation*}
        J(t_0, x_0, u(\cdot))
        \geq J_\eta(t_0, z_0, p_\eta(\cdot)) - \varepsilon / 3
        \geq \Valaeta(t_0, z_0) - \varepsilon / 3.
    \end{equation*}
    Since this estimate is valid for every $u(\cdot) \in \U(t_0, \vartheta)$, we get $\Val(t_0, x_0) \geq \Valaeta(t_0, z_0) - \varepsilon / 3$.

    On the other hand, let $p_\eta(\cdot) \in \U(t_0, \vartheta_\eta)$, and let $u(t) = p_\eta(t)$ for $t \in [t_0, \vartheta_\eta)$ and $u(t) = \bar{u} \in \USet$ for $t \in [\vartheta_\eta, \vartheta)$ (see (\ref{Theorem_OCP_varepsilon_main})).
    Arguing as above, we obtain
    \begin{equation} \label{Theorem_OCP_varepsilon_proof_main_2}
        J_\eta(t_0, z_0, p_\eta(\cdot))
        \geq J(t_0, x_0, u(\cdot)) - \varepsilon / 3
        \geq \Val(t_0, x_0) - \varepsilon /3.
    \end{equation}
    Consequently, $\Valaeta(t_0, z_0) \geq \Val(t_0, x_0) - \varepsilon / 3$, and, hence,
    \begin{equation} \label{Theorem_OCP_varepsilon_proof_Val_Valeta}
        |\Val(t_0, x_0) - \Valaeta(t_0, z_0)|
        \leq \varepsilon / 3
        \leq \varepsilon.
    \end{equation}

    Now, let a control $p^\ast_\eta(\cdot)$ be $\varepsilon^\ast$-optimal in the auxiliary problem, and let the control $u^\ast(\cdot)$ be defined by (\ref{Theorem_OCP_varepsilon_main}).
    Then, in accordance with (\ref{Theorem_OCP_varepsilon_proof_main_2}) and (\ref{Theorem_OCP_varepsilon_proof_Val_Valeta}), we have
    \begin{equation*}
        J(t_0, x_0, u^\ast(\cdot))
        \leq J_\eta(t_0, z_0, p^\ast_\eta(\cdot)) + \varepsilon / 3
        \leq \Valaeta(t_0, z_0) + \varepsilon^\ast + \varepsilon / 3
        \leq \Val(t_0, x_0) + \varepsilon.
    \end{equation*}
    Thus, the control $u^\ast(\cdot)$ is $\varepsilon$-optimal in the original problem.
    The theorem is proved.
\qed

{\it Proof of Theorem~\ref{Theorem_OFC}}
    Let $\varepsilon > 0$ be fixed.
    Let us choose $\eta^\ast \in (0, \vartheta - t_0)$ as in the proof of Theorem~\ref{Theorem_OCP_varepsilon} and put $\eta^\circ = \min\{\eta^\ast, \varepsilon / 3\}$.
    For every $\eta \in (0, \eta^\circ]$, let us determine $\omega^\circ_\eta(\eta, \varkappa^\circ_\eta(\eta))$ according to (\ref{kappa_delta_eta_a}) and define $\delta^\circ(\eta) = \omega^\circ_\eta(\eta, \varkappa^\circ_\eta(\eta))$.

    Let us fix $\eta \in (0, \eta^\circ]$, a partition $\Delta$ (\ref{Delta}) such that $\diam(\Delta) \leq \delta^\circ(\eta)$, and $x_0 \in B(R_x)$, and consider the control $u(\cdot) = u(\cdot \mid t_0, x_0, \vartheta, U^\ast, \eta, \Delta)$ formed in system (\ref{system}) by the control law $\{U^\ast, \eta, \Delta\}$ on the basis of the control strategy $U^\ast$ from (\ref{strategy_U_ast}).
    Thus, in order to complete the proof, it is sufficient to show that (see (\ref{optimal_control_strategy}))
    \begin{equation} \label{Theorem_OFC_proof_final}
        J(t_0, x_0, u(\cdot))
        \leq \Val(t_0, x_0) + \varepsilon.
    \end{equation}

    Let us consider the partition $\Delta_\eta = \{\tau_j\}_{j = 1}^{k + 1}$ of $[t_0, \vartheta_\eta]$ (see (\ref{Delta_eta})) such that
    \begin{equation*}
        \Delta_\eta = (\Delta \cap [t_0, \vartheta_\eta]) \cup \{\vartheta_\eta\}.
    \end{equation*}
    Let us note that $\diam(\Delta_\eta) \leq \diam(\Delta) \leq \delta^\circ(\eta) = \omega^\circ_\eta(\eta, \varkappa^\circ_\eta(\eta))$.
    Let $z_0 = \Image(t_0, z_0)$, and let the control $p_\eta(\cdot) = p_\eta(\cdot \mid t_0, z_0, \vartheta_\eta, P_\eta^\circ, \varkappa^\circ_\eta(\eta), \Delta_\eta)$ be formed in auxiliary system (\ref{system_z_eta}) by the control law $\{P_\eta^\circ, \varkappa^\circ_\eta(\eta), \Delta_\eta\}$ on the basis of the optimal in the auxiliary problem control strategy $P_\eta^\circ$.
    Then, in accordance with (\ref{I_x_0_R_z}) and (\ref{optimal_control_strategy_a}), taking (\ref{Theorem_OCP_varepsilon_proof_Val_Valeta}) into account, we obtain
    \begin{equation} \label{Theorem_OFC_proof_1}
        J_\eta (t_0, z_0, p_\eta(\cdot))
        \leq \Valaeta(t_0, z_0) + \eta
        \leq \Valaeta(t_0, z_0) + \varepsilon / 3
        \leq \Val(t_0, x_0) + 2 \varepsilon / 3.
    \end{equation}

    Further, arguing by induction, let us prove that
    \begin{equation} \label{proof_induction}
        u(t)
        = p_\eta(t),
        \quad t \in [\tau_j, \tau_{j + 1}), \quad j \in \overline{1, k}.
    \end{equation}
    For $j = 1$, due to (\ref{Control_law}), (\ref{Control_law_eta}), and (\ref{strategy_U_ast}), we have
    \begin{equation*}
        u(t)
        = U^\ast(t_0, x_0, \eta)
        = P^\circ_\eta \big(t_0, z_0, \varkappa^\circ_\eta (\eta) \big)
        = p_\eta(t),
        \quad t \in [\tau_1, \tau_2).
    \end{equation*}
    Now, let us take $q \in \overline{2, k}$ and suppose that (\ref{proof_induction}) holds for every $j \in \overline{1, q - 1}$, i.e., $u(t) = p_\eta(t)$, $t \in [t_0, \tau_{q})$.
    Hence, $\Image(\tau_q, x_{\tau_q}(\cdot)) = z(\tau_q)$ by Lemma~\ref{Lemma_Connection}, and, therefore, we get (\ref{proof_induction}) for $j = q$:
    \begin{equation*}
        u(t)
        = U^\ast(\tau_q, x_{\tau_q}(\cdot), \eta)
        = P^\circ_\eta \big(\tau_q, z(\tau_q), \varkappa^\circ_\eta (\eta) \big)
        = p_\eta(t),
        \quad t \in [\tau_q, \tau_{q + 1}).
    \end{equation*}

    Applying (\ref{proof_induction}) for $j = k$, we conclude $u(t) = p_\eta(t)$, $t \in [t_0, \vartheta_\eta)$.
    Then, in accordance with (\ref{Theorem_OCP_varepsilon_proof_main_2}) and (\ref{Theorem_OFC_proof_1}), we derive
    \begin{equation*}
        J(t_0, x_0, u(\cdot))
        \leq J_\eta(t_0, z_0, p_\eta(\cdot)) + \varepsilon / 3
        \leq \Val(t_0, x_0) + \varepsilon.
    \end{equation*}
    Thus, the inequality in (\ref{Theorem_OFC_proof_final}) and the theorem are proved.
\qed

\bibliographystyle{spmpsci_unsrt}
\bibliography{C:/Users/Mikhail/Dropbox/TeX/BibTeX/references}

\end{document}